\documentclass[12pt,a4paper,twoside]{article}
\usepackage{amsmath,amsfonts,amssymb,a4,enumerate,epsfig,array,theorem}

\newtheorem{theo}{Theorem}
\newtheorem{prop}{Proposition}[section]
\newtheorem{coro}[prop]{Corollary}

\newtheorem{lemma}[prop]{Lemma}

\newcommand{\ZZ}{{\mathbb{Z}}}
\newcommand{\QQ}{{\mathbb{Q}}}
\newcommand{\RR}{{\mathbb{R}}}

\newcommand{\tD}{{\tilde D}}

\newcommand\rank{\operatorname{rank}}
\newcommand\step{\operatorname{step}}

\newcommand{\nobf}{\noindent\bf}

\def\qed{\unskip\nobreak\hfil\penalty50\hskip1.75em\null\nobreak\hfil
$\blacksquare$ {\parfillskip=0pt \finalhyphendemerits=0 \par}\goodbreak}

\begin{document}
\title{Arithmetic properties of the adjacency matrix \\ of quadriculated disks}
\author{Nicolau C. Saldanha and Carlos Tomei}
\maketitle

\begin{abstract}
Let $\Delta$ be a bicolored quadriculated disk
with black-to-white matrix $B_\Delta$.
We show how to factor $B_\Delta = L\tilde DU$,
where $L$ and $U$ are lower and upper triangular matrices,
$\tilde D$ is obtained from a larger identity matrix
by removing rows and columns and
all entries of $L$, $\tilde D$ and $U$ are
equal to $0$, $1$ or $-1$.
\end{abstract}

\section{Introduction}

\begin{figure}[ht]
\begin{center}
\epsfig{height=3cm,file=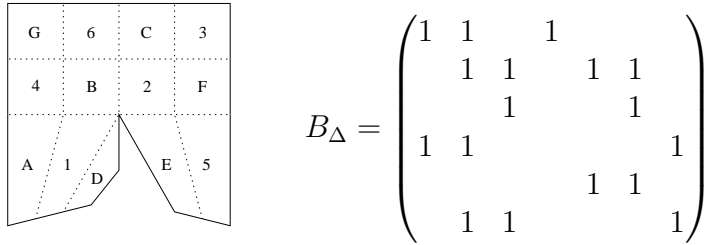}
\qquad
\raise 12mm \hbox{\(
B_\Delta = \begin{pmatrix}
1 & 1 &  & 1 &  &  &  \\
 & 1 & 1 &  & 1 & 1 &  \\
 &  & 1 &  &  & 1 &  \\
1 & 1 &  &  &  &  & 1 \\
 &  &  &  & 1 & 1 &  \\
 & 1 & 1 &  &  &  & 1
\end{pmatrix}
\)}
\end{center}
\label{fig:disk0}
\caption{A quadriculated disk and its black-to-white matrix}
\end{figure}

\footnote{1991 {\em Mathematics Subject Classification}.
Primary 05B45, Secondary 05A15, 05C50, 05E05.
{\em Keywords and phrases} Quadriculated disk, tilings by dominoes,
dimers.}
\footnote{The authors gratefully acknowledge the support of
CNPq, Faperj and Finep.}
In this paper, a \textit{square} is a topological disk
with four privileged boundary points, the \textit{vertices}.
A \textit{quadriculated disk} is a closed topological disk
formed by the juxtaposition along sides of finitely many squares
such that interior vertices belong to precisely four squares.
A simple example of a quadriculated disk is the $n \times m$ rectangle
divided into unit squares, another example is given in figure \ref{fig:disk0}.
A quadriculated disk may be considered as a closed subset of the plane $\RR^2$.

Quadriculated disks are {\it bicolored},
i.e., the squares are black and white in a way that
squares with a common side have opposite colors.
Label the black (resp. white) squares of a quadriculated disk $\Delta$
by $1, 2, \ldots, b$ (resp. $1, 2, \ldots, w$).
The $b \times w$ {\it black-to-white} matrix $B_\Delta$
has $(i,j)$ entry $b_{ij} = 1$ if the $i$-th black and $j$-th white squares
share an side and $b_{ij} = 0$ otherwise.
If the disk has a single square of either color the matrix $B$ collapses.
See figure \ref{fig:disk0} for an example of black-to-white matrix;
black squares are numbered.
Throughout the paper, blank matrix entries equal $0$.
For a given labeling of squares of a quadriculated disk $\Delta$
in which black squares come first, the adjacency matrix is
\[ \begin{pmatrix} 0 & B \\ B^T & 0 \end{pmatrix}. \]

A rectangular matrix $\tD$ is a \textit{defective identity}
if it can be obtained from the identity matrix
by adding rows and columns of zeroes or,
equivalently, by removing rows and columns from a larger identity matrix.
More precisely, $\tD$ is a defective identity if every entry equals $0$ or $1$
and if $\tilde d_{ij} = \tilde d_{i'j'} = 1$ implies either
$i = i'$, $j = j'$, or $i < i'$, $j < j'$, or $i > i'$, $j > j'$.
For a $n \times m$ matrix $A$,
an $L\tD U$ decomposition of $A$ is a factorization $A = L\tD U$
where $L$ (resp. $U$) is $n \times n$ (resp. $m \times m$)
lower (resp. upper) invertible and $\tD$ is a defective identity.
The main result of this paper is the following.

\begin{theo}
\label{theo:LDU}
Let $\Delta$ be a quadriculated disk with at least two squares.
For an appropriate labeling of its squares,
the black-to-white matrix $B_\Delta$
admits an $L\tD U$ decomposition
whose factors have all entries equal to $0$, $1$ or $-1$.
\end{theo}

Thus, for example, the matrix $B_\Delta$ in figure \ref{fig:disk0}
admits the decomposition
\[
\begin{pmatrix}
1 & & & & & \\
& 1 & & & & \\
& & 1 & & & \\
1 & & & -1 & & \\
& & & & 1 & \\
& 1 & & & -1 & 1 
\end{pmatrix}
\begin{pmatrix}
1 & & & & & & \\
& 1 & & & & & \\
& & 1 & & & & \\
& & & 1 & & & \\
& & & & 1 & 0 & \\
& & & & & & 1
\end{pmatrix}
\begin{pmatrix}
1 & 1 & & 1 & & & \\
& 1 & 1 & & 1 & 1 & \\
& & 1 & & & 1 & \\
& & & 1 & & & -1 \\
& & & & 1 & 1 & \\
& & & & & 1 & \\
& & & & & & 1
\end{pmatrix}.
\]

The determinant $\det(B_\Delta)$ is, up to sign, the $q$-counting
for $q = -1$ of the domino tilings of $\Delta$, in the sense of \cite{EKLP}.
A {\it domino tiling} of $\Delta$ is a decomposition of $\Delta$ as a disjoint
union of $2 \times 1$ rectangles: $b=w$ is a necessary but not sufficient
condition for the existence of a domino tiling.
The result below was originally proved in \cite{DT}
using a very different argument.

\begin{coro}
\label{coro:DT}
Let $\Delta$ be a quadriculated disk with $b=w$ and black-to-white
matrix $B_\Delta$. Then $\det(B_\Delta)$ equals $0$, $1$ or $-1$.
\end{coro}

Theorem \ref{theo:LDU} is stronger, however,
especially when $B_\Delta$ is not invertible.

\begin{coro}
\label{coro:tor}
Let $\Delta$ be a quadriculated disk with black-to-white matrix $B_\Delta$.
If $v$ has integer entries and the system $B_\Delta x = v$ admits
a rational solution then the system admits an integer solution.
\end{coro}

This corollary may be interpreted as saying that the cokernel
$\ZZ^b/B_\Delta(\ZZ^w)$ of $B_\Delta: \ZZ^w \to \ZZ^b$ is a free abelian group.
From Theorem \ref{theo:LDU},
the rank $r$ of $B_\Delta$ is the same in $\QQ$ as in $\ZZ_p$
for any prime number $p$.
The proof of Corollary \ref{coro:DT} in \cite{DT}
is based on this fact for $p=2$.

The proof of Theorem \ref{theo:LDU} is constructive in the sense
that it yields a fast algorithm to obtain the appropriate labeling
of vertices, the matrices in the factorization,
$\det(B_\Delta)$ and $\rank(B_\Delta)$.

In section \ref{sect:diags} we introduce \textit{cut and paste},
a geometric procedure that converts a quadriculated disk $\Delta$ into
a union of smaller quadriculated disks $\Delta'_1,\ldots,\Delta'_\ell$.
In section \ref{sect:decomp} we present its algebraic counterpart,
relating the black-to-white matrices $B_{\Delta'_1},\ldots,B_{\Delta'_\ell}$
to the original $B_\Delta$; this will be the inductive step in the proof
of Theorem \ref{theo:LDU}.
In section \ref{sect:boards} we study \textit{boards},
quadriculated disks which are subsets of the quadriculated plane
$\ZZ^2 \subset \RR^2$, and show that cut and paste can be performed within
this smaller class (Theorem \ref{theo:board}).

We would like to thank Tania V. de Vasconcelos for helpful conversations
and the referee for a careful report which led to a clearer
and more correct paper.

\section{Cut and paste}
\label{sect:diags}

\begin{figure}[ht]
\begin{center}
\epsfig{height=35mm,file=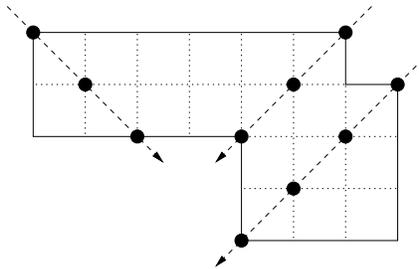}
\end{center}
\label{fig:gooddiag}
\caption{Good balanced, bad and good unbalanced diagonals}
\end{figure}

The proof of Theorem \ref{theo:LDU} makes use of \textit{good diagonals}:
examples are given in figure \ref{fig:gooddiag}.
\textit{Cut and paste} along a good diagonal
obtains from a quadriculated disk $\Delta$
a disjoint union of smaller disks $\Delta'_1, \ldots, \Delta'_\ell$,
often with $\ell = 1$ (Lemma \ref{lemma:deltalinha}).
Our argument consists of two steps: a proof that there always exists
a good diagonal (Proposition \ref{prop:gooddiag}) and
a procedure to convert $L\tilde DU$ decompositions of
$B_{\Delta'_1}, \ldots, B_{\Delta'_\ell}$
into a decomposition of $B_\Delta$ (Lemma \ref{lemma:LDUstep}).

Topological subdisks of $\RR^2$
consisting of unit squares with vertices in $\ZZ^2$ are {\it boards}:
the quadriculated disk in figure \ref{fig:disk0} is not a board.
Arbitrary quadriculated disks can still be ``ironed''
so that its squares become unit squares with vertices in $\ZZ^2$:
this is what \textit{developing maps} make precise
(see \cite{TL} for a more general discussion of developing maps).

\begin{lemma}
\label{lemma:devel}
Let $\Delta$ be a quadriculated disk:
there exists a continuous map $\phi: \Delta \to \RR^2$
which is locally injective in the interior of $\Delta$ and
takes squares in $\Delta$ to unit squares in $\RR^2$ with vertices in $\ZZ^2$.
\end{lemma}

{\nobf Proof: }
Orient $\Delta$ and endow it with a metric structure such that each square
becomes isometric to the unit square in the plane:
$\phi$ is to be an orientation preserving local isometry
in the interior of $\Delta$.
Pick an arbitrary edge $v_0v_1$ and set $\phi(v_0) = (0,0)$,
$\phi(v_1) = (1,0)$ and extend $\phi$ isometrically along the edge.
Notice that an orientation preserving local isometry
extends uniquely from an edge to the squares adjacent to the edge.
By extending the domain of $\phi$ repeatedly starting from $v_0v_1$,
uniqueness is clear.
Since $\Delta$ is simply connected, $\phi$ is well-defined.
\qed

For the quadriculated disk in figure \ref{fig:disk0},
the developing map $\phi$ is neither injective in the interior of $\Delta$
nor locally injective on its boundary.

There exist genuinely different quadriculated disks $\Delta$
and $\tilde \Delta$ for which the image of the restriction
of the developing map $\phi$ to their boundaries is the same:
this a corollary of Milnor's paisley curve example (\cite{P}).
It is not hard to produce such examples with
$\det(B_\Delta) = 0 \neq \det(B_{\tilde \Delta})$.

A \textit{corner} of a quadriculated disk
is a boundary point which is a vertex of a single square.
A \textit{diagonal} of a quadriculated disk is
a sequence of vertices $v_0v_1\ldots v_k$, $k > 0$, such that
({\it i}) {$v_0$ is a corner, $v_1, v_2, \ldots, v_{k-1}$ are interior vertices
and $v_k$ is a boundary vertex;}
({\it ii}) {consecutive vertices $v_i$ and $v_{i+1}$, $i = 0, \ldots, k-1$,
are opposite vertices of a square $s_{i+1}$;}
({\it iii}) {consecutive squares $s_i$ and $s_{i+1}$, $i = 1, \ldots, k-1$,
have a single vertex in common (which is $v_i$).}

Notice that both coordinates of $\phi(v_i)$ depend monotonically on $i$.
In particular,
vertices and squares of a diagonal of a quadriculated disk are all distinct.
A diagonal $v_0\ldots v_k$ is a \textit{good diagonal}
if at least one side of the square $s_k$
containing $v_k$ is contained in the boundary of $\Delta$.
A good diagonal $v_0\ldots v_k$ is \textit{balanced}
(see figure \ref{fig:gooddiag})
if exactly one side of $s_k$ containing $v_k$
is contained in the boundary of $\Delta$,
and \textit{unbalanced} otherwise.

\begin{prop}
\label{prop:gooddiag}
Any quadriculated disk $\Delta$ admits at least four good diagonals.
\end{prop}

Notice that diagonals, being sequences of vertices,
are naturally oriented. In particular,
a square has four diagonals, all good.

{\nobf Proof:}
Let $V$, $E$ and $F$ be the number of vertices, edges and faces,
i.e., squares, of $\Delta$.
Write $E = E_I + E_B$, where $E_I$ (resp. $E_B$)
counts interior (resp. boundary) edges:
clearly, $4F = 2E_I + E_B = 2E - E_B$
and therefore $4E = 8F + 2E_B$.
Also, $V = V_I + V_1 + V_2 + \cdots + V_r$
where $V_r$ is the number of boundary vertices
belonging to $r$ squares.
Again, $4F = 4 V_I + V_1 + 2 V_2 + \cdots + r V_r =
4V - (3 V_1 + 2 V_2 + \cdots + (4-r) V_r)$
and
$4V = 4F + (3V_1 + 2V_2 + \cdots + (4-r) V_r)$.

By Euler, $4V - 4E + 4F = 4$.
Substituting the above formulae and using
$E_B = V_1 + V_2 + \cdots + V_r$ we have
$V_1 - 4 = V_3 + 2V_4 + \cdots + (r-2)V_r$.

Each vertex counted in $V_1$ is a starting corner for a diagonal:
we have to prove that at least four of these $V_1$ diagonals are good.
Each vertex counted in $V_3$, for example, is the endpoint
of at most three diagonals of which only one is declared bad.
More generally, we have at most $V_1 - 4 = V_3 + 2V_4 + \cdots + (r-2)V_r$
bad ends and we are done.
\qed

\begin{figure}[ht]
\begin{center}
\epsfig{height=45mm,file=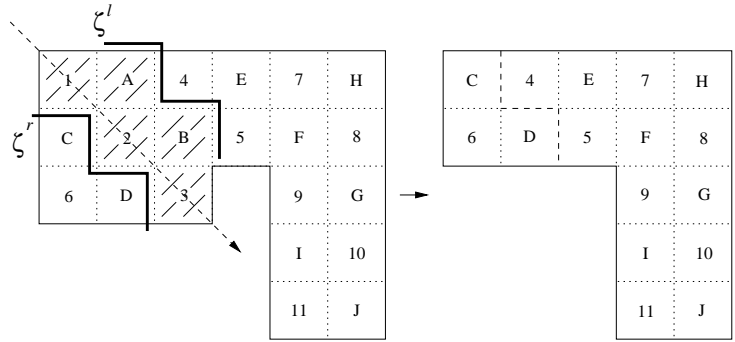}
\end{center}
\caption{Cut and paste along an unbalanced diagonal}
\label{fig:cutnpaste0}
\end{figure}

\begin{figure}[ht]
\begin{center}
\epsfig{height=55mm,file=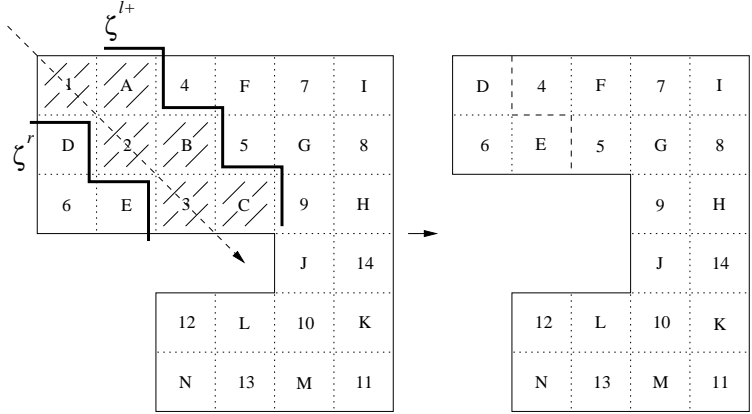}
\end{center}
\caption{Cut and paste along a balanced diagonal}
\label{fig:cutnpaste1}
\end{figure}

We are ready to {\it cut and paste} along a good diagonal.
Figures \ref{fig:cutnpaste0} and \ref{fig:cutnpaste1}
illustrate cut and paste along balanced and unbalanced diagonals.
Let $\zeta^l$, $\zeta^{l+}$ and $\zeta^r$
be the zig-zags indicated in the figures.
In figure \ref{fig:cutnpaste0},
cut and paste removes the shaded squares
and identifies $\zeta^l$ and $\zeta^r$
to obtain a new quadriculated disk.
In figure \ref{fig:cutnpaste1},
$\zeta^r$ is identified with the first four segments of $\zeta^{l+}$
(near the top).

\begin{figure}[ht]
\begin{center}
\epsfig{height=4cm,file=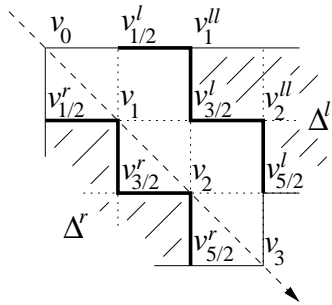}
\end{center}
\caption{Notation for cut and paste}
\label{fig:indices}
\end{figure}

More precisely, cut and paste obtains from
a quadriculated square $\Delta$ a disjoint union of quadriculated squares
$\Delta' = \Delta'_1 \cup \cdots \cup \Delta'_\ell$.
In the two previous examples, $\ell = 1$;
also, as we shall see, $\Delta' = \emptyset$ if and only if
$\Delta$ consists of one or two squares.
Consider a good diagonal with squares $s_1, \ldots, s_k$.
If the diagonal is balanced, choose one side of the diagonal
so as to contain $k$ squares; otherwise choose any side.
For convenience, we always assume the chosen side to be the left,
the other case being analogous. 
Squares to the left of the diagonal will be called $s^l_1, \ldots$;
the last square in this sequence is $s^l_{k'}$:
$k' = k$ if the diagonal is balanced and $k' = k-1$ otherwise.
Squares to the right of the diagonal will be called
$s^r_1, \ldots, s^r_{k-1}$;
notice that there is no $s^r_k$.
For $k > 1$, label the vertices of $s_i$ counterclockwise
$v_{i-1}$, $v^r_{i-1/2}$, $v_i$, $v^l_{i-1/2}$
and the vertices of $s^l_i$ counterclockwise
$v^l_{i-1/2}$, $v_i$, $v^l_{i+1/2}$, $v^{ll}_i$
(see figure \ref{fig:indices}).
The zig-zag $\zeta^l$ is
$v^l_{1/2}, v^{ll}_1, v^l_{3/2}, \ldots, v^{ll}_{k-1}, v^l_{k-1/2}$
and $\zeta^r$ is
$v^r_{1/2}, v_1, v^r_{3/2}, \ldots, v_{k-1}, v^r_{k-1/2}$.
If the diagonal is unbalanced then $\zeta^{l+} = \zeta^l$;
otherwise $\zeta^{l+}$ is two segments longer and ends with
$v^{ll}_k, v^l_{k+1/2}$.
Both $\zeta^{l+}$ and $\zeta^r$ begin and end at boundary points.
The region between $\zeta^{l+}$ and $\zeta^r$ consists precisely
of the squares $s_1,\ldots,s_k$ and $s^l_1,\ldots,s^l_{k'}$.
Let $\Delta^r$ (resp. $\Delta^l$) be the closed regions
to the right (resp. left) of $\zeta^r$ (resp. $\zeta^l$).
Attach $\Delta^l$ to $\Delta^r$
by identifying $\zeta^l$ and $\zeta^r$:
in order to obtain a quadriculated region $\tilde\Delta'$.

\begin{lemma}
\label{lemma:deltalinha}
For $k > 1$, the region $\Delta^r$ constructed above
is non-empty, path connected and simply connected;
$\Delta^l$ is a (possibly empty) disjoint union
of path connected and simply connected regions.
The region $\tilde\Delta'$ is connected and simply connected
and therefore a union of quadriculated disks
$\Delta'_1, \ldots, \Delta'_\ell$ with $\Delta'_i \cap \Delta'_j$, $i \ne j$,
consisting of at most a single vertex.
If $k = 1$, $\Delta'$ is obtained from $\Delta$
by excising the two squares $s_1$ and $s^l_1$.
\end{lemma}

\begin{figure}[ht]
\begin{center}
\epsfig{height=55mm,file=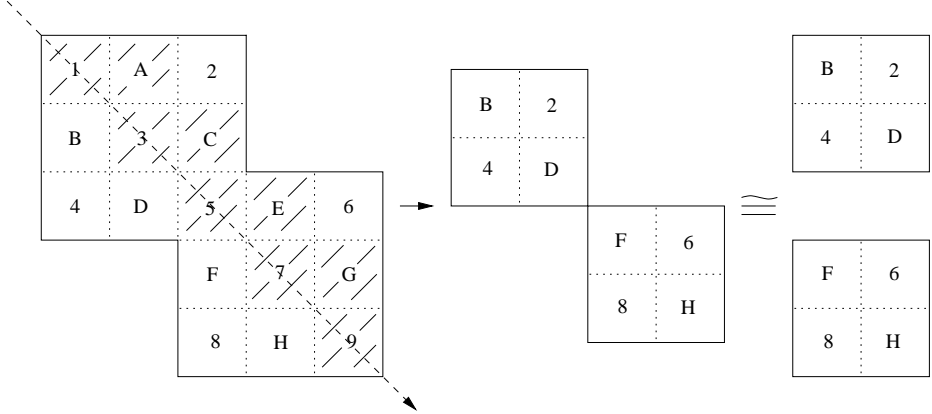}
\end{center}
\caption{Cut and paste may produce a disjoint union of disks}
\label{fig:cutnpaste2}
\end{figure}

{\nobf Proof:}
Notice that
some isolated points in $\zeta^r$ of the form $v^r_i$ may be boundary
points but otherwise $\zeta^r$ is in the interior of $\Delta$;
the zig-zag $\zeta^{l+}$, on the other hand, may contain boundary segments
(see figure \ref{fig:cutnpaste2}).
To prove path connectivity of $\Delta^r$,
take the shortest path in $\Delta$ connecting an arbitrary point
in $\Delta^r$ to a point in $\zeta^r$.
If a simple closed curve is contained in $\Delta^r$ (or $\Delta^l$)
then so is the disk surrounded by the curve:
this prove simple connectivity of the connected components.
Since $\Delta^l$, $\Delta^r$ and their intersection $\zeta^l = \zeta^r$
in $\tilde\Delta'$ are all path connected and simply connected, so is $\Delta'$.
Interior vertices of $\tilde\Delta'$ are surrounded by four squares 
and therefore $\tilde\Delta'$ is quadriculated.
\qed

With the notation of the lemma,  set $\Delta'$
to be the disjoint union of the quadriculated disks
$\Delta'_1,\ldots, \Delta'_\ell$.
In this sense, cut and paste obtains from a disk 
a (possibly empty) disjoint union of quadriculated disks.
Under the obvious extension of the concept
of black-to-white matrix, 
\[ B_{\tilde\Delta'} = B_{\Delta'}
= \begin{pmatrix} B_{\Delta'_1} & \ldots & 0 \\
\vdots & & \vdots \\ 0 & \ldots & B_{\Delta'_\ell} \end{pmatrix} \]
for the same labeling of squares
and thus $L\tilde DU$ decompositions of
$B_{\Delta'_1}, \ldots, B_{\Delta'_\ell}$
yields a similar decomposition of $B_{\Delta'}$.
Figure \ref{fig:cutnpaste2} shows an example of $\Delta$,
$\tilde\Delta'$ and $\Delta'$ where $\ell = 2$.

\vfill\eject

\section{The decomposition}
\label{sect:decomp}

Before discussing the relationship between $B_\Delta$ and $B_{\Delta'}$
we present a lemma in linear algebra.
The proof is a straightforward computation left to the reader.
Let $I_n$ be the $n \times n$ identity matrix and
$I_{n,m}$ be the $n \times m$ defective identity matrix
with $(i,j)$ entry equal to 1 if $i=j$ and 0 otherwise.

\begin{lemma}
\label{lemma:linalg}
Let $M$ be a $(n + m) \times (n' + m')$ matrix which can be
decomposed in blocks as
\[ M = \begin{pmatrix} M_{11} & M_{12} \\ M_{21} & M_{22} \end{pmatrix}, \]
where $M_{11}$ is $n \times n'$.
If $n' \le n$ and $N$ is a $n' \times m'$ matrix with $M_{11} N = M_{12}$ then
\[ M =
\begin{pmatrix} M_{11} I_{n',n} & 0 \\ M_{21} I_{n',n} & I_{m} \end{pmatrix}
\begin{pmatrix} I_{n,n'} & 0 \\ 0 & M_{22} - M_{21} N \end{pmatrix}
\begin{pmatrix} I_{n'} & N \\ 0 & I_{m'} \end{pmatrix}.
\]
Similarly, if $n' \ge n$ and $N$ is a $m \times n$ matrix
with $N M_{11} = M_{21}$ then
\[ M =
\begin{pmatrix} I_{n} & 0 \\ N & I_{m} \end{pmatrix}
\begin{pmatrix} I_{n,n'} & 0 \\ 0 & M_{22} - N M_{12} \end{pmatrix}
\begin{pmatrix} I_{n',n} M_{11} &
I_{n',n} M_{12} \\ 0 & I_{m'} \end{pmatrix}. 
\]
\end{lemma}

The next lemma is the inductive step in the proof of Theorem \ref{theo:LDU}.

\begin{lemma}
\label{lemma:LDUstep}
Let $\Delta$ be a quadriculated disk with $b$ black and $w$ white squares,
$b+w > 1$.
Let $\Delta' = \Delta'_1 \cup \cdots \cup \Delta'_\ell$
(with $b' = b'_1 + \cdots + b'_\ell$ black
and $w' = w'_1 + \cdots + w'_\ell$ white squares)
be obtained from $\Delta$ by cut and paste along a good diagonal.
Label black and white squares in $\Delta$ so that excised squares
come first and in the order prescribed by the good diagonal;
label squares in $\Delta'$ next.
Then the black-to-white matrices $B_\Delta$ and $B_{\Delta'}$ satisfy
\[ B_\Delta =
\begin{pmatrix} L & 0 \\ X & S_{b'} \end{pmatrix}
\begin{pmatrix} I_{b-b',w-w'} & 0 \\ 0 & B_{\Delta'} \end{pmatrix}
\begin{pmatrix} U & Y \\ 0 & S_{w'} \end{pmatrix}
\]
where $L$ (resp $U$) is an invertible lower (resp. upper)
square matrix of order $b-b'$ (resp. $w-w'$)
and $S_{b'}$ and $S_{w'}$ are square diagonal matrices.
Furthermore, all entries of $S_{b'}$, $S_{w'}$, $L$, $U$, $X$ and $Y$
equal 0, 1 or $-1$.
\end{lemma}

The statement above requires clarification in some degenerate cases.
If $\Delta'$ is empty, $B_{\Delta'}$ collapses and
$B_\Delta = L I_{b,w} U$.
If instead $\Delta'$ is a disjoint union of unit squares,
all of the same color, then either $w' = 0$ or $b' = 0$ and
\[ B_\Delta =
\begin{pmatrix} L & 0 \\ X & S_{b'} \end{pmatrix}
\begin{pmatrix} I_{b-b',w} \\ 0 \end{pmatrix} U
\quad\hbox{or}\quad
B_\Delta =
L \begin{pmatrix} I_{b,w-w'} & 0 \end{pmatrix}
\begin{pmatrix} U & Y \\ 0 & S_{w'} \end{pmatrix}.
\]

{\nobf Proof:}
Assume that the squares on the good diagonal are black; thus $k = b-b'$;
if the squares were white all computations would be transposed.
Let $j_1, \ldots, j_{k-1}$ be the indices of the white squares
$s^r_1, \ldots, s^r_{k-1}$; notice that $j_i > w-w'$.
Decompose the matrix $B_\Delta$ in four blocks,
\[ B_\Delta = \begin{pmatrix} B_{11} & B_{12} \\ B_{21} & B_{22} \end{pmatrix}, \]
where $B_{22}$ is a $b' \times w'$ matrix.
By construction, $B_{11}$ has one of the two forms below,
the first case corresponding to balanced good diagonals
(i.e., to $b-b' = w-w'$).
\[
B_{11} = \begin{pmatrix}
1 & 0 & 0 & \cdots & 0 & 0 \\
1 & 1 & 0 & \cdots & 0 & 0 \\
0 & 1 & 1 & \cdots & 0 & 0 \\
\vdots & \vdots & \vdots & & \vdots & \vdots \\
0 & 0 & 0 & \cdots & 1 & 0 \\
0 & 0 & 0 & \cdots & 1 & 1
\end{pmatrix},
\quad \hbox{or} \quad
B_{11} = \begin{pmatrix}
1 & 0 & 0 & \cdots & 0  \\
1 & 1 & 0 & \cdots & 0  \\
0 & 1 & 1 & \cdots & 0  \\
\vdots & \vdots & \vdots & & \vdots \\
0 & 0 & 0 & \cdots & 1  \\
0 & 0 & 0 & \cdots & 1 
\end{pmatrix}.
\]
Let $S_b$ (resp. $S_w$) be a $b \times b$ (resp. $w \times w$)
diagonal matrix with diagonal entries equal to $1$ or $-1$;
the $i$-th entry of $S_b$ (resp. $S_w$) is $-1$ if the $i$-th
black (resp. white) square is strictly to the right of the diagonal.
Write
\[ S_b = \begin{pmatrix} I_{b-b'} & 0 \\ 0 & S_{b'} \end{pmatrix}, \qquad
S_w = \begin{pmatrix} I_{w-w'} & 0 \\ 0 & S_{w'} \end{pmatrix}. \]
We have
\[ S_b B_\Delta S_w =
\begin{pmatrix} B_{11} & - B_{12} \\ B_{21} & B_{22} \end{pmatrix}. \]
The matrix $B_{12}$ has $k-1$ nonzero columns
in positions $j_1, \ldots, j_{k-1}$.
More precisely, the nonzero entries are
$(i,j_i)$ and $(i+1,j_i)$ for $i = 1,\ldots,k-1$.
The nonzero columns of $B_{12}$ are equal to the first $k-1$ columns
of $B_{11}$. Let $N$ be the $(w-w') \times w'$ matrix with entries 0 or $-1$,
with nonzero entries at $(1,j_1), (2,j_2), \ldots, (k-1,j_{k-1})$.
Clearly $B_{11} N = - B_{12}$ and we may apply Lemma \ref{lemma:linalg}
to write $S_b B_\Delta S_w$ as
\[
\begin{pmatrix} B_{11} I_{w-w',b-b'} & 0 \\
B_{21} I_{w-w',b-b'} & I_{b'} \end{pmatrix}
\begin{pmatrix} I_{b-b',w-w'} & 0 \\ 0 & B_{22} - B_{21} N \end{pmatrix}
\begin{pmatrix} I_{w-w'} & N \\ 0 & I_{w'} \end{pmatrix}.
\]
We claim that $B_{\Delta'} = B_{22} - B_{21} N$.
The nonzero columns of the matrix $- B_{21} N$ are
the columns of $B_{21}$, except that the first column is moved
to position $j_1$, the second column is moved to $j_2$ and so on.
The $k$-th column of $B_{21}$, if it exists, is discarded.
These nonzero entries correspond precisely to the identifications
which must be performed in order to obtain $\Delta'$, i.e.,
to the ones which must be added to $B_{22}$ in order to obtain $B_{\Delta'}$.
We now have
\[ B_\Delta =
\begin{pmatrix} B_{11} I_{w-w',b-b'} & 0 \\
S_{b'} B_{21} I_{w-w',b-b'} & S_{b'} \end{pmatrix}
\begin{pmatrix} I_{b-b',w-w'} & 0 \\ 0 & B_{\Delta'} \end{pmatrix}
\begin{pmatrix} I_{w-w'} & N S_{w'} \\ 0 & S_{w'} \end{pmatrix}.
\]
If the good diagonal is balanced, this finishes the proof.
In the unbalanced case, $\tilde L = B_{11} I_{w-w',b-b'}$ is not invertible
since its last column is zero. Replace the $(k,k)$ entry of 
$\tilde L$ by 1 to obtain a new matrix $L$:
$L$ is clearly invertible and $\tilde L I_{b-b',w-w'} = L I_{b-b',w-w'}$.
The proof is now complete.
\qed

{\nobf Proof of Theorem \ref{theo:LDU}:}
The basis of the induction on the number of squares of $\Delta$
consists of checking that the theorem holds for disks with at most two squares.
Notice that if the disk consists of a single square then $b = 0$ or $w = 0$
and the matrices are degenerate.

Let $\Delta$ be a quadriculated disk and
$\Delta' = \Delta'_1 \cup \cdots \cup \Delta'_\ell$ 
be obtained from $\Delta$ by cut and paste.
By induction on the number of squares the theorem may be assumed to hold
for eack $\Delta'_k$ and we therefore write
$B_{\Delta'} = L_{\Delta'} \tD_{\Delta'} U_{\Delta'}$.
From the induction step, Lemma \ref{lemma:LDUstep}, write
\begin{align}
B_\Delta &=
\begin{pmatrix} L_{\step} & 0 \\ X_{\step} & S_{b'} \end{pmatrix}
\begin{pmatrix} I_{b-b',w-w'} & 0 \\ 0 & B_{\Delta'} \end{pmatrix}
\begin{pmatrix} U_{\step} & Y_{\step} \\ 0 & S_{w'} \end{pmatrix}
\notag\\
&= L_\Delta \tD_\Delta U_\Delta. \notag
\end{align}
where
\[
L_\Delta =
\begin{pmatrix} L_{\step} & 0 \\ X_{\step} & S_{b'} \end{pmatrix}
\begin{pmatrix} I_{b-b'} & 0 \\ 0 & L_{\Delta'} \end{pmatrix},
\quad
U_\Delta =
\begin{pmatrix} I_{w-w'} & 0 \\ 0 & U_{\Delta'} \end{pmatrix}
\begin{pmatrix} U_{\step} & Y_{\step} \\ 0 & S_{w'} \end{pmatrix}.
\]
The theorem now follows from observing that each nonzero entry of $L_\Delta$
(resp. $U_\Delta$)
is, up to sign, copied from either
$L_{\Delta'}$, $L_{\step}$ or $X_{\step}$
(resp. $U_{\Delta'}$, $U_{\step}$ or $Y_{\step}$)
and is therefore equal to $1$ or $-1$.
\qed

\begin{figure}[ht]
\begin{center}
\epsfig{height=50mm,file=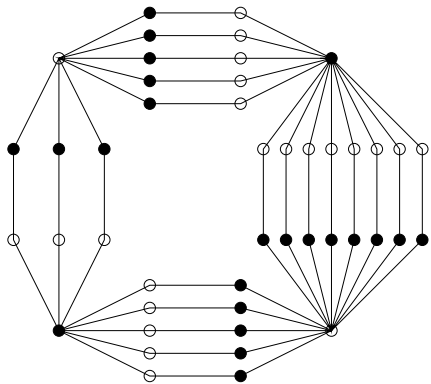}
\end{center}
\caption{Determinant 1 does not imply $L\tilde DU$ decomposition}
\label{fig:fibo}
\end{figure}

The example in Figure \ref{fig:fibo} is instructive:
the $B_G$ matrix of this planar graph $G$ has determinant 1
but admits no $L\tilde DU$ decomposition where the matrices
have integer coefficients since the removal of any two vertices
of opposite colors from $G$ yields a graph whose determinant
has absolute value greater than 1.

\section{Boards}
\label{sect:boards}

\begin{figure}[ht]
\begin{center}
\epsfig{height=4cm,file=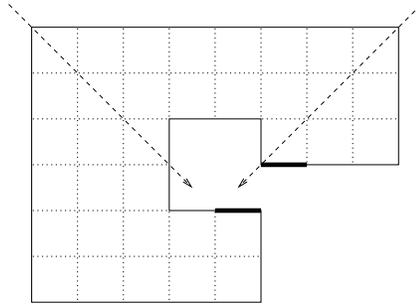}
\end{center}
\caption{A board and two good diagonals, one excellent.}
\label{fig:board}
\end{figure}

A \textit{board} is a quadriculated disk $\Delta$ for which
the developing map $\phi: \Delta \to \RR^2$ is injective.
In other words, a board is a topological subdisk of $\RR^2$
whose boundary is a polygonal curve consisting of segments
of length 1 joining points in $\ZZ^2$.
In this section we always assume $\phi$ to be the inclusion.
The class of boards is not closed under cut and paste:
in figure \ref{fig:board}, the two enhanced segments on the boundary
would be superimposed by cut and paste with the good diagonal on the left.
Cut and paste along the good diagonal indicated on the right, however,
yields a smaller board.

\begin{theo}
\label{theo:board}
It is always possible to cut and paste a given board $\Delta$
to obtain a disjoint union of boards $\Delta'$.
\end{theo}

{\nobf Proof:}
A diagonal is \textit{excellent} if the $x$ and $y$ coordinates
are both monotonic along one of the two boundary arcs between $v_0$ and $v_k$;
without loss, let this arc lie to the right of the diagonal.
Excellent diagonals are good: the vertex $v^r+{k-1/2}$ is on the boundary.
We interpret cut and paste along an excellent diagonal
as leaving $\Delta^l$ fixed and moving $\Delta^r$.
In this way, $\Delta'$ becomes a subset of $\Delta$
and is therefore a disjoint union of boards.
We are left with proving that any board admits excellent diagonals.

Each diagonal defines two boundary arcs:
order these arcs by inclusion.
We claim that a diagonal defining a minimal arc is excellent.
Indeed, let $\delta^m = (v^m_0v^m_1\ldots v^m_k)$
be a diagonal inducing a minimal arc $\alpha$:
assume without loss of generality that $v^m_i = (a+i,b+i)$ 
for given $a$ and $b$.
Assume by contradiction that $\delta^m$ is not excellent.
Consider the value of $x+y$ on $\alpha$:
this function is not strictly increasing and therefore admits a local
minimum somewhere in the interior of $\alpha$.
The local minimum closest to $v^m_0$ is the starting point for a diagonal
defining a boundary arc $\beta$ strictly contained in $\alpha$,
a contradiction.
\qed

\vfill

\goodbreak

\bigskip\bigskip\bigbreak

{

\parindent=0pt
\parskip=0pt
\obeylines

Nicolau C. Saldanha and Carlos Tomei 
Departamento de Matem\'atica, PUC-Rio
R. Marqu\^es de S. Vicente 225, Rio de Janeiro, RJ 22453-900, Brazil

\medskip

nicolau@mat.puc-rio.br; http://www.mat.puc-rio.br/$\sim$nicolau/
tomei@mat.puc-rio.br

}

\end{document}